\documentclass[11pt,twoside]{sat}
\newcommand{\sect}[1]{\section{#1}\setcounter{equation}{0}}

\title{Pad\'e and Hermite-Pad\'e approximation and orthogonality}
\def\shorttitle{Pad\'e, Hermite-Pad\'e, and orthogonality}

\author{Walter Van Assche}
\def\shortauthor{Walter Van Assche}

\def\versiondate{1 September 2006}

\def\startpagenumber{61}
\def\volumenumber{2}
\def\year{2006}
\newcommand{\beginddoc}{
\begin{document}
\maketitle
\begin{abstract}
We give a short introduction to Pad\'e approximation (rational approximation to a function
with close contact at one point) and to Hermite-Pad\'e approximation (simultaneous rational
approximation to several functions with close contact at one point) and show how orthogonality
plays a crucial role. We give some insight into how logarithmic potential theory helps in describing
the asymptotic behavior and the convergence properties of Pad\'e and Hermite-Pad\'e approximation. 

MSC: 41A21, 41A28, 42C05, 33C45
\end{abstract}
\insert\footins{\scriptsize
\medskip
\baselineskip 8pt
\leftline{Surveys in Approximation Theory}
\leftline{Volume \volumenumber, \year.
pp.~\thepage--\pageref{endpage}.}
\leftline{Copyright \copyright\ 2005 Surveys in Approximation Theory.}
\leftline{ISSN 1555-578X}
\leftline{All rights of reproduction in any form reserved.}
\smallskip
\par\allowbreak}
\tableofcontents}
\markboth{{\it \shortauthor}}{{\it \shorttitle}}
\markright{{\it \shorttitle}}
\def\endddoc{\label{endpage}\end{document}}
\date{{\small \versiondate}}
\def\dword#1{{\bf#1}}

\usepackage{amsfonts}
\usepackage{amsmath}
\usepackage{amsthm}
\renewcommand{\O}{\mathcal{O}}
\renewcommand{\L}{\mathcal{L}}
\newtheorem{theorem}{Theorem}[section]
\newtheorem{definition}{Definition}[section]
\newtheorem{lemma}{Lemma}[section]

\input colordvi
\def\comment#1{{\bf(comment: #1)}}
\def\comment#1{\Blue{{\bf \  #1}}}
\def\marginal#1{\llap{#1\ }}
\def\ct#1#2{#2}
\def\ct#1#2{{\bf[}#1{\bf]\hbox{\tt -->}[}#2{\bf]}}
\def\cte#1#2{\Green{{\bf[}\Red{#1}{\bf]}\hbox{$\;\Rightarrow\;$}{\bf[}\Red{#2}{\bf]}}}
\def\ct#1#2{\cte{#1}{#2}\marginpar{!}}

\beginddoc


\sect{Pad\'e approximation}
\subsection{Taylor polynomials}
The general setup in approximation theory is that a function $f$
is given and that one wants to approximate it with a
\textit{simpler} function $g$ but in such a way that the
difference between $f$ and $g$ is small. The advantage is that the
simpler function $g$ can be handled without too many difficulties
but the disadvantage is that one loses some information since $f$
and $g$ are different.

In the setting of Pad\'e approximation one starts with a
function $f: \mathbb{C} \to \mathbb{C}$ for which a Taylor
expansion is known in the neighborhood of a given point $a \in
\mathbb{C}$, i.e.,
\begin{equation}  \label{eq:fexp}
f(z) = \sum_{k=0}^\infty c_k (z-a)^k , \quad c_k = \frac{f^{(k)}(a)}{k!}.
\end{equation}
The function $f$ can not be computed exactly using this Taylor
expansion since this requires an infinite number of additions (and
multiplications). We can obtain a polynomial approximation
by truncating after $n$ terms. The corresponding approximations
are Taylor polynomials given by
\begin{equation}
f_n(z) = \sum_{k=0}^{n-1} c_k (z-a)^k, 
\end{equation}
and these Taylor polynomials are therefore characterized by
\begin{equation}
    f(z) - f_n(z) = \O((z-a)^n), \qquad z \to a.
\end{equation}
This condition is a (confluent) interpolation condition which
tells us that the difference $f-f_n$ has a zero of multiplicity
$n$ at the point $a$. We know an explicit formula for the Taylor
polynomial, namely
\[   f_n(z) = \sum_{k=0}^{n-1} \frac{f^{(k)}(a)}{k!} (z-a)^k, \]
and the error is given by
\[   f(z)-f_n(z) = \sum_{k=n}^\infty \frac{f^{(k)}(a)}{k!}
(z-a)^k. \] If $f$ is analytic in a domain $\Omega$ that contains
$a$ and if $\Gamma$ is a closed contour in $\Omega$ encircling $a$
once in the positive direction (counterclockwise), then Cauchy's
formula gives
\[     \frac{f^{(k)}(a)}{k!} = \frac{1}{2\pi i} \int_{\Gamma}
\frac{f(\xi)}{(\xi-a)^{k+1}} \, d\xi, \]
and hence
\begin{eqnarray*}
   f_n(z) & = & \frac{1}{2\pi i} \int_{\Gamma} \frac{f(\xi)}{\xi-a}
  \sum_{k=0}^{n-1} \left( \frac{z-a}{\xi-a} \right)^k \ d\xi \\
     & = & \frac{1}{2\pi i} \int_{\Gamma} \frac{f(\xi)}{\xi-z} \left[ 1-
     \left( \frac{z-a}{\xi-a} \right)^n \right] \ d\xi.
 \end{eqnarray*}
 The error then becomes
 \begin{equation}  \label{eq:Taylorer}
    f(z)-f_n(z) = \frac{1}{2\pi i} \int_{\Gamma} \frac{f(\xi)}{\xi-z}  \left( \frac{z-a}{\xi-a} \right)^n  \
 d\xi.
 \end{equation}
 The convergence of $f_n$ to $f$ corresponds to the convergence of
 the Taylor series, and typically one has uniform convergence on
 closed disks $|z-a| \leq r$, where $r < \rho(f)$ and
 \[    \rho(f) := \sup \{ R : f \textrm{ is analytic in } |z-a|<R \}  \]
is the radius of convergence of the series in (\ref{eq:fexp}). 
Indeed, if we choose
$\epsilon > 0$ such that $r+\epsilon < \rho(f)$ and if we take for
$\Gamma$ the circle $|\xi-a|=r+\epsilon$, then for $|z-a|\leq r$
we have from (\ref{eq:Taylorer}) by straightforward estimations
\[  |f(z)-f_n(z)| \leq  \max_{|\xi| = r+\epsilon}
|f(\xi)| \left( \frac{r}{r+\epsilon} \right)^{n} \frac{1}{2\pi}
\int_\Gamma \frac{|d\xi|}{|\xi-z|} , \]
 and since $r/(r+\epsilon) < 1$ we see that the right hand side
converges to $0$. So convergence is only guaranteed on disks with
a radius less than the radius of convergence. The function $f$ may be
analytic in a larger domain (the radius of convergence depends on the
singularity of $f$ closest to $a$), but the Taylor approximation
will not converge outside the disk with radius $\rho(f)$.

\subsection{Pad\'e approximants}
Polynomials are not such a good class of functions if one wants to
approximate functions with singularities because polynomials are
entire functions without singularities. They are only useful up to
the first singularity of $f$ near $a$. Rational functions are the
simplest functions with singularities. The idea is that the poles
of the rational functions will move to the singularities of the
function $f$, and hence the domain of convergence could be
enlarged, and singularities of $f$ may be discovered using the
poles of the rational approximants.

The $[m,n]$ Pad\'e approximant of $f$ in $a$ is the rational
function $Q_m/P_n$, with $Q_m$ a polynomial of degree $\leq m$ and
$P_n$ a polynomial of degree $\leq n$, for which we have the
following interpolation condition at $a$:
\begin{equation}
  f(z) - \frac{Q_m(z)}{P_n(z)} = \O((z-a)^{m+n+1}), \qquad z\to a.
\end{equation}
The computation of the polynomials $P_n$ and $Q_m$ is not so easy
from this interpolation condition, since one first has the compute
the Taylor expansion of $Q_m/P_n$ and then equate the first
$m+n+1$ Taylor coefficients to the first $m+n+1$ Taylor
coefficients of $f$. Usually the Pad\'e approximant is defined by
linearizing the interpolation condition as
\begin{equation}  \label{eq:pade}
    P_n(z)f(z) - Q_m(z) = \O((z-a)^{m+n+1}), \qquad z \to a.
\end{equation}
For Pad\'e approximation near infinity to a function of the form
\[   f(z) = \sum_{k=0}^\infty \frac{c_k}{z^{k+1}}, \]
one takes $m=n-1$ and the interpolation condition is
\[   P_n(z) f(z) - Q_{n-1}(z) = \O(z^{-n-1}), \qquad z \to \infty, \]
(see Section 1.3).
There is a degree of freedom since we can multiply both sides of
(\ref{eq:pade}) by a constant. Usually we normalize this by taking
$P_n$ \dword{monic} (i.e., of the form $x^n+\cdots$) when this is
possible, and this can only be done if $P_n$ is of exact degree
$n$.  If we take $P_n$ monic, then we can determine the $n$
unknown coefficients $a_k$ $(k=1,\ldots,n)$ in 
\begin{equation}  \label{eq:Pnexp}
   P_n(z) =: \sum_{k=0}^n a_k (z-a)^{n-k}, \qquad a_0=1,
\end{equation} by putting the coefficients of $(z-a)^k$ for
$k=m+1,m+2,\ldots,m+n$ in the Taylor expansion of $P_nf$
equal to zero. The polynomial $Q_m$ then corresponds to the Taylor
polynomial of degree $m$ of $P_nf$.

Here is another approach. Suppose $f$ is analytic in a domain
$\Omega$ that contains $a$. Again we take a contour $\Gamma$
inside $\Omega$ encircling $a$ once in the positive direction.
Divide both sides of (\ref{eq:pade}) by $(z-a)^{m+k+2}$ and
integrate, to find
\begin{multline*}
   \frac{1}{2\pi i} \int_\Gamma \frac{P_n(z)f(z)}{(z-a)^{m+k+2}}
\, dz  - \frac{1}{2\pi i} \int_\Gamma \frac{Q_m(z)}{(z-a)^{m+k+2}}
\,dz \\
= \sum_{j=m+n+1}^\infty  b_{n,j} \frac{1}{2\pi i} \int_\Gamma
(z-a)^{j-m-k-2}\, dz,
\end{multline*} 
where the $b_{n,j}$'s are the coefficients in the expansion of 
$P_nf-Q_m$ around $a$.
The integral involving $Q_m$ is
zero for $k\ge0$ since it is proportional to the $(m+k+1)$th derivative of
$Q_m$, which is zero for $k \geq 0$. The sum on the right-hand
side has a contribution only when $j=m+k+1$, but when $0\leq k
\leq n-1$ then $j \leq m+n$ and such indices do not appear in the
sum. Hence the right hand side also vanishes for $k\le n-1$. Therefore
(\ref{eq:pade}) implies that
\[  \frac{1}{2\pi i} \int_\Gamma \frac{P_n(z)}{(z-a)^{m+k+2}}
f(z)\, dz = 0, \qquad k=0,1,\ldots,n-1. \]
If we use the expansion
(\ref{eq:Pnexp}) then this gives
\[  \sum_{j=0}^n a_j \frac{1}{2\pi i} \int_\Gamma (z-a)^{n-j-m-k-2}
f(z) \, dz = 0, \qquad k=0,1,2,\ldots,n-1. \] If we use the
expansion (\ref{eq:fexp}) then
\[  \frac{1}{2\pi i} \int_\Gamma (z-a)^{n-j-m-k-2}
f(z) \, dz = c_{m-n+k+j+1}, \] so we get the system of equations
\begin{equation}  \label{eq:system}
\begin{pmatrix}
        c_{m-n+1} & c_{m-n+2} & \cdots & c_{m+1} \\
        c_{m-n+2} & c_{m-n+3} & \cdots & c_{m+2} \\
        \vdots & \vdots & \cdots & \vdots \\
        c_{m} & c_{m+1} & \cdots & c_{m+n}
       \end{pmatrix}
       \begin{pmatrix} a_0 \\ a_1 \\ \vdots \\ a_{n}
       \end{pmatrix} =
       \begin{pmatrix} 0 \\ 0 \\ \vdots \\ 0
      \end{pmatrix} .
\end{equation}
There is one degree of freedom here since we have $n+1$ unknowns
and $n$ (homogeneous) equations. The choice $a_0=1$ (if possible) gives the
monic polynomial $P_n$, but sometimes another normalization will
be used, as we will see later.

\subsection{Orthogonality}
 From now on we will only consider Pad\'e approximants near
infinity. This can easily be obtained from Pad\'e approximation near zero
and the change of variable $z \mapsto 1/z$. Indeed, if $g$ has a Taylor expansion
\[    f^*(z) := \sum_{k=0}^\infty c_k z^k   \]
near the origin, then $f(z):=g(1/z)/z$ as an expansion near infinity of the form
\begin{equation}  \label{eq:finf}
    f(z)= \sum_{k=0}^\infty  \frac{c_k}{z^{k+1}}.
\end{equation}
Since $f(z) = \O(1/z)$, the only sensible choice of the
degree in the rational approximation problem is to take $m=n-1$ so
that $Q_m(z)/P_n(z)$ is also $\O(1/z)$. This situation occurs when $f$ is
of the form
\[   f(z) = \int_{-\infty}^\infty \frac{d\mu(x)}{z-x}, \]
i.e., when $f$ is the Stieltjes transform (or Cauchy transform) of a positive measure $\mu$ on
the real line. The Pad\'e approximants
near infinity can be obtained from the Pad\'e approximants near zero
in the following way.  The $[n-1,n]$ Pad\'e approximant
$Q^*_{n-1}/P^*_n$ for $f^*$ near 0 has the interpolation condition
\[   P^*_n(x) f^*(x) - Q^*_{n-1}(x) = \O(x^{2n}), \qquad x \to 0.
\]
Change variables by setting $x=1/z$ and divide both sides by 
$z$. Then
\[   P^*_n(1/z) f(z) - \frac1z Q^*_{n-1}(1/z) = \O(z^{-2n-1}),
\qquad z \to \infty. \] In order to get polynomials, we multiply
both sides by $z^n$. Then
\begin{equation}  \label{eq:padeinf}
   P_n(z) f(z) - Q_{n-1}(z) = \O(z^{-n-1}) , \qquad z \to \infty,
\end{equation}
 where $P_n(z) := z^n P_n^*(1/z)$ and
$Q_{n-1}(z) := z^{n-1} Q_{n-1}^*(1/z)$ are obtained by reversing
the polynomials $P_n^*$ and $Q_{n-1}^*$. So the interpolation
conditions at infinity are given by (\ref{eq:padeinf}). The system of equations
(\ref{eq:system}) for $f^*$ and $m=n-1$ then changes to the system
\begin{equation}  \label{eq:system2}
\begin{pmatrix}
        c_{0} & c_{1} & \cdots & c_{n} \\
        c_{1} & c_{2} & \cdots & c_{n+1} \\
        \vdots & \vdots & \cdots & \vdots \\
        c_{n} & c_{n+1} & \cdots & c_{2n-1}
       \end{pmatrix}
       \begin{pmatrix} a_0 \\ a_{1} \\ \vdots \\ a_{n}
       \end{pmatrix} =
       \begin{pmatrix} 0 \\ 0 \\ \vdots \\ 0
      \end{pmatrix} ,
\end{equation}
for the unknown coefficients of
\[  P_n(z) := \sum_{k=0}^n a_k z^k.  \]

Typically we will not be given the function $f$ but rather the
infinite sequence of coefficients $c_0,c_1,c_2,\ldots$ in the
Laurent expansion of $f$. With this as input, we define a linear
functional $\L$ on the linear space of polynomials by
\begin{equation}  \label{eq:L}
   \L(x^n) := c_n, \qquad n=0,1,2,\ldots.
\end{equation}
For a polynomial $p(x)=\sum_{k=0}^n a_k x^k$ we then have by
linearity $\L(p) = \sum_{k=0}^n a_kc_k$. If we now look at the
system of equations (\ref{eq:system2}), then the
coefficients of $P_n$ satisfy the
equations
\[       \sum_{j=0}^n a_j c_{k+j} = 0, \qquad k=0,1,\ldots,n-1. \]
But this is equivalent to saying that
\begin{equation}  \label{eq:Lortho}
   \L(x^k P_n(x)) = 0, \qquad k=0,1,\ldots,n-1.
\end{equation}
Hence the polynomial $P_n$ is orthogonal to all
polynomials of degree less than $n$ with respect to the linear functional
$\L$. A very useful normalization of $P_n$ is to require 
that in addition to (\ref{eq:Lortho}) we also have
\[   \L(P_n^2(x)) = 1. \]
This can always be done when the functional is positive. 
When the functional is not positive, then one imposes the extra
condition $\L(P_n^2(x)):= h_n \neq 0$, so that $P_n/\sqrt{h_n}$ has norm one. 
Once the polynomial $P_n$
is obtained, the remaining elements in the Pad\'e approximation
problem can be found explicitly in terms of $P_n$. Indeed, if we define 
\begin{equation}  \label{eq:Q}
    Q_{n-1}(z) := \L \left( \frac{P_n(z)-P_n(x)}{z-x} \right),
 \end{equation}
then, since  $[P_n(z)-P_n(x)]/(z-x)$ is a polynomial of degree
$n-1$ in the variable $z$, $Q_{n-1}$ is a polynomial of degree $n-1$
and (\ref{eq:Q}) is equivalent to
\[    P_n(z) \L \left( \frac{1}{z-x} \right) - Q_{n-1}(z) = \L
\left( \frac{P_n(x)}{z-x} \right). \]
The functional $\L$ was only defined on polynomials, but if we expand
$1/(z-x)$ in a Laurent series, then (at least
formally)
\[     \L \left( \frac{1}{z-x} \right) = \L \left(
\sum_{k=0}^\infty \frac{x^k}{z^{k+1}} \right) = \sum_{k=0}^\infty
\frac{c_k}{z^{k+1}} = f(z), \] so what needs to be shown is that
\[   \L \left( \frac{P_n(x)}{z-x} \right) = \O(z^{-n-1}). \]
Using the Laurent series of $1/(z-x)$ we find
\[    \L \left( \frac{P_n(x)}{z-x} \right) = \sum_{k=0}^\infty
\frac{1}{z^{k+1}} \L(x^k P_n(x)), \] 
and the orthogonality
conditions (\ref{eq:Lortho}) show that the terms with $k \leq n-1$
vanish. The first term is therefore the term with $k=n$, which is
$\O(1/z^{n+1})$. What we also learn from this proof is that the error
in the Pad\'e approximation problem is given explicitly by
\begin{equation}  \label{eq:error}
     P_n(z) f(z) - Q_{n-1}(z) = \L \left( \frac{P_n(x)}{z-x}
     \right),
\end{equation}
which is again in terms of the polynomial $P_n$.

\subsection{Moment problem}
The linear functional $\L$ remains a bit mysterious. Obviously it is
related to the function $f$, but we would like to know it somewhat
more explicitly. The Riesz representation theorem tells us that
every positive and bounded linear functional on the linear space
of continuous functions with compact support
 on the real line can be represented by a
finite positive measure $\mu$ on the real line as
\[  \L(f) = \int_{-\infty}^\infty f(x)\, d\mu(x). \]
If we want to get convergence results for Pad\'e approximation,
then it would be convenient to work with a bounded and positive
linear functional $\L$, which is represented by a finite positive
measure $\mu$. In that case
\begin{equation} \label{eq:moments}
   c_k =\int_{-\infty}^\infty x^k\, d\mu(x)
\end{equation}
 will be the moments of a positive measure $\mu$ and the
function $f$ is the Cauchy transform (Stieltjes transform) of the measure $\mu$:
\[    f(z) = \int_{-\infty}^\infty \frac{1}{z-x} \, d\mu(x). \]
Obviously not every infinite sequence $c_0,c_1,c_2,\ldots$ will
lead to a positive and bounded linear functional. The moment
problem is to obtain conditions on this infinite sequence
$c_0,c_1,c_2,\ldots$ guaranteeing
 that they are the moments of a
finite positive measure on the real line, as in
(\ref{eq:moments}). If the measure is supported on
$(-\infty,\infty)$ then this is known as the \dword{Hamburger
moment problem}. If the measure is supported on the positive axis
$[0,\infty)$ then we speak of the \dword{Stieltjes moment
problem}. If the measure is supported on a finite interval
(usually $[0,1]$), then this is known as the \dword{Hausdorff
moment problem}. A necessary and sufficient condition that the
sequence $c_0,c_1,c_2,\ldots$ consist of moments of a positive
measure on $(-\infty,\infty)$ is that all the \dword{Hankel
matrices}
\[  \begin{pmatrix}
        c_{0} & c_{1} & \cdots & c_{n} \\
        c_{1} & c_{2} & \cdots & c_{n+1} \\
        \vdots & \vdots & \cdots & \vdots \\
        c_{n} & c_{n+1} & \cdots & c_{2n}
       \end{pmatrix} \]
be positive definite. Observe that these are precisely the
matrices appearing in (\ref{eq:system2}).

 From now on we will add one more restriction, namely that the
measure be supported on a finite interval $[a,b]$. This
simplifies our treatment by avoiding non-compactness of the
support. So our function $f$ will be a \dword{Markov function}
\[   f(z) = \int_a^b \frac{1}{z-x} \, d\mu(x), \]
and such a function is analytic in $\mathbb{C} \setminus [a,b]$.
The singularities of this function therefore are located on the
interval $[a,b]$. The linear functional in this case is given by
\[   \L(g) = \int_a^b g(x)\, d\mu(x), \]
for every continuous function $g$ on $[a,b]$.
The denominator polynomials in the Pad\'e
approximation problem are orthogonal polynomials for the measure
$\mu$ on the interval $[a,b]$, i.e.,
\begin{equation} \label{eq:muortho}
  \int_a^b  x^k P_n(x)\, d\mu(x) = 0, \qquad k=0,1,\ldots,n-1,
  \end{equation}
  which we normalize so that they are orthonormal
\begin{equation}  \label{eq:norm1}
  \int_a^b P_n^2(x)\, d\mu(x) = 1.
 \end{equation}
 The numerator polynomials are given by
\begin{equation}  \label{eq:Qnum}
   Q_{n-1}(z) = \int_a^b \frac{P_n(z)-P_n(x)}{z-x}\, d\mu(x),
  \end{equation}
 and the error is given by
 \begin{equation}  \label{eq:errormu}
     P_n(z) f(z) - Q_{n-1}(z) = \int_a^b \frac{P_n(x)}{z-x} \,
     d\mu(x).
  \end{equation}

\subsection{Zeros and poles}
The idea of using rational approximation is that the singularities
of the Pad\'e approximant would give an idea of the singularities
of the function $f$. This is indeed so when $f$ is a Markov
function. The singularities of the Pad\'e approximant are  poles
at the zeros of $P_n$. A consequence of the orthogonality is that
these zeros are simple and they all are on the open interval
$(a,b)$.

\begin{theorem}
Suppose that the support of $\mu$ is an infinite set in $[a,b]$.
Then all the zeros of $P_n$ are simple and located on $(a,b)$.
\end{theorem}

\begin{proof}
Let $x_1,\ldots,x_m$ be the sign changes of $P_n$ on $(a,b)$, then
obviously $m \leq n$, since each sign change is a zero. Suppose
that $m < n$. Then introduce the polynomial
$\pi_m(x) := (x-x_1)(x-x_2)\cdots(x-x_m)$. The function
$P_n(x)\pi_m(x)$ does not change sign on $[a,b]$ and since the
support of $\mu$ contains infinitely many points we conclude that
\[ \int_a^b P_n(x)\pi_m(x)\, d\mu(x) \neq 0. \]
But $P_n$ is orthogonal to all polynomials of degree $<n$,
hence this integral is equal to 0. This contradiction implies that
$m=n$. So $P_n$ has $n$ sign changes on $(a,b)$, each a zero of
$P_n$, hence each a simple zero of $P_n$, and $P_n$ has no other zeros.
\end{proof}

\subsection{Convergence}
When we study the convergence of the Pad\'e approximants, we use
(\ref{eq:errormu}) to find
\[   f(z)- \frac{Q_{n-1}(z)}{P_n(z)} = \frac{1}{P_n(z)} \int_a^b
\frac{P_n(x)}{z-x} \, d\mu(x). \] 
Observe that
\begin{multline*}
   P_n(z) \int_a^b \frac{P_n(x)}{z-x} \, d\mu(x) \\ 
= \int_a^b
\frac{P_n(x)[P_n(z)-P_n(x)]}{z-x}\, d\mu(x) + \int_a^b
\frac{P_n^2(x)}{z-x}\, d\mu(x) . 
\end{multline*}
 The fraction
$[P_n(z)-P_n(x)]/(z-x)$ is a polynomial of degree $n-1$ in the
variable $x$, so by orthogonality the first integral on the right
vanishes. This gives
\[  P_n(z) \int_a^b \frac{P_n(x)}{z-x} \, d\mu(x) =  \int_a^b \frac{P_n^2(x)}{z-x} \,
d\mu(x), \] and the error in Pad\'e approximation becomes
\begin{equation}  \label{eq:errorP2}
  f(z)- \frac{Q_{n-1}(z)}{P_n(z)} = \frac{1}{P_n^2(z)} \int_a^b
\frac{P_n^2(x)}{z-x} \, d\mu(x).
\end{equation}
This error contains two parts: on the one hand it contains the
polynomial $P_n$ for which we will describe the asymptotic
behavior in the next subsection, and on the other hand it contains
the integral
\[   \int_a^b \frac{P_n^2(x)}{z-x} \, d\mu(x), \]
which is in fact a Markov function for the probability measure
$P_n^2(x)\, d\mu(x)$ when $P_n$ is the orthonormal polynomial. We
can estimate this integral as follows. Suppose that $z$ belongs to
a compact set $K \subset \mathbb{C} \setminus [a,b]$. Then the
distance $d_K$ between $K$ and $[a,b]$
\[    d_K = \inf \{ |z-x| : z\in K, x \in [a,b] \} \]
is strictly positive. Therefore we have
\[  \left|  \int_a^b \frac{P_n^2(x)}{z-x}\, d\mu(x) \right| \leq
\int_a^b \frac{P_n^2(x)}{|z-x|} \, d\mu(x) \leq \frac{1}{d_K},
\]
and this bound is independent of $n$. So the convergence of the
Pad\'e approximants is completely determined by the asymptotic
behavior of $P_n$.

\subsection{Asymptotic properties}
In this subsection we describe the asymptotic behavior of
$|P_n(z)|^{1/n}$ when $z \in K$, where $K$ is a compact subset of
$\mathbb{C} \setminus [a,b]$. If we denote the leading coefficient
of $P_n$ by $\gamma_n >0$ and the zeros of $P_n$ by $x_{1,n} <
x_{2,n} < \cdots < x_{n,n}$, then
\[   P_n(z) = \gamma_n \prod_{j=1}^n (z-x_{j,n}). \]
The asymptotic behavior thus requires knowing the behavior of
$\gamma_n$ and the asymptotic distribution of the zeros.

Let us first consider the asymptotic distribution of the zeros.
Consider the discrete measure
\[  \nu_n := \frac1n \sum_{j=1}^n \delta_{x_{j,n}}, \]
where $\delta_c$ is the Dirac measure with mass $1$ at the point
$c$. The measure $\nu_n$ describes the distribution of the zeros
of $P_n$. The asymptotic distribution corresponds to an
investigation of the limit of this sequence of measures. All the
zeros of $P_n$ are on the interval $[a,b]$, so all the measures
$\nu_n$ are probability measures on $[a,b]$. Helly's selection
principle tells us that there will be a subsequence that converges
weakly to a probability measure $\nu$ on $[a,b]$. This means that
there is a subsequence $(n_k)$ such that
\[  \lim_{k \to \infty} \int_a^b g(x)\, d\nu_{n_k}(x) = \int_a^b
g(x)\, d\nu(x), \] for every continuous function $g$ on $[a,b]$.
For the monic polynomial $\hat{P}_n := P_n/\gamma_n$ we have
\[  \frac1n \log |\hat{P}_n(z)| = \frac1n \sum_{j=1}^n \log
|z-x_{j,n}| = \int_a^b \log |z-x|\, d\nu_n(x), \] hence when $z
\in K \subset \mathbb{C}\setminus [a,b]$, then the weak
convergence implies that
\[  \lim_{k \to \infty} |\hat{P}_{n_k}(z)|^{1/n_k} = \exp \left(
\int_a^b \log|z-x|\, d\nu(x) \right). \]

Next, the leading coefficient $\gamma_n$ solves a minimization
problem:
\begin{theorem}
We have
 \begin{equation}  \label{eq:gamma}
   \frac{1}{\gamma_n^2} = \min_{q_n(x)=x^n+\cdots} \int_a^b
   |q_n(x)|^2\, d\mu(x),
 \end{equation}
and the minimum is attained at the monic orthogonal polynomial
$\hat{P}_n$.
\end{theorem}

\begin{proof}
We can write an arbitrary monic polynomial of degree $n$ as $q_n =
\hat{P}_n + \pi_{n-1}$, where $ \pi_{n-1}$ is a polynomial
of degree $\leq n-1$. We then have
\begin{multline*}
  \int_a^b |q_n(x)|^2\, d\mu(x) = \int_a^b |\hat{P}_n(x)|^2\,
d\mu(x) + \int_a^b |\pi_{n-1}(x)|^2\, d\mu(x) \\
 + 2 \int_a^b
\hat{P}_n(x)\pi_{n-1}(x) \, d\mu(x). 
\end{multline*}
The last integral vanishes
because of  orthogonality, so that
\[  \min_{q_n(x)=x^n+\cdots}\int_a^b |q_n(x)|^2\, d\mu(x)
    = \int_a^b |\hat{P}_n(x)|^2\, d\mu(x) + \min_{\pi_{n-1}}
    \int_a^b |\pi_{n-1}(x)|^2\, d\mu(x). \]
The minimum on the right hand side is obtained by
 taking  $\pi_{n-1}=0$, so the minimum in (\ref{eq:gamma}) is
 obtained for the monic orthogonal polynomial.
 \end{proof}

Without going to much into details, this extremal problem for
$\gamma_n$ will in fact tell us that the asymptotic behavior of
$\gamma_n^{1/n}$ and the asymptotic distribution of the zeros (the
measure $\nu$) are described by an equilibrium problem for
(logarithmic) potentials. There is a unique probability measure
$\mu_e$ on $[a,b]$  that minimizes the \dword{logarithmic energy}
\[  \int_a^b \int_a^b \log \frac{1}{|x-y|} \,d\sigma(x)d\sigma(y)
\] over all probability measures $\sigma$ supported on $[a,b]$.
This measure is given by
\[  d\mu_e(x) = \frac{1}{\pi} \frac{dx}{\sqrt{(x-a)(b-x)}}, \quad
x \in [a,b] \] and has the property that its \dword{logarithmic potential} satisfies
\[  U(x;\mu_e) = \int_a^b \log \frac{1}{|x-y|} \, d\mu_e(y) = -\log
\frac{b-a}4, \qquad x \in [a,b]. \] This \dword{equilibrium
measure} corresponds to the measure $\nu$ describing the
asymptotic zero distribution when the orthogonality measure $\mu$
is sufficiently regular on $[a,b]$. A sufficient condition is that
$\mu' > 0$ almost everywhere on $[a,b]$ (\dword{Erd\H{o}s-Tur\'an
condition}). Furthermore, we also have
\[  \lim_{n \to \infty} \gamma_n^{1/n} = \frac4{b-a}. \]
Combining both results shows that when $\mu'>0$ almost everywhere
on $[a,b]$ we have
\[  \lim_{n \to \infty} |P_n(z)|^{1/n} = \frac{4}{b-a} \exp \left(
  -\int_a^b \log \frac{1}{|z-x|} \, d\mu_e(x) \right). \]
When $z$ is on the interval $[a,b]$ then the right hand side is
equal to $1$, but when $z$ moves away from $[a,b]$, then the right
hand side becomes $>1$. On the equipotential curves
\[  C_r = \{ z \in \mathbb{C} \setminus [a,b]: \frac{4}{b-a} \exp \left(- \int_a^b \log
\frac{1}{|z-x|}\, d\mu_e(x) \right) = r \} \] with $r> 1$ we then
conclude that
\[  \lim_{n \to \infty} | f(z) - \frac{Q_{n-1}(z)}{P_n(z)}|^{1/n} =
\frac{1}{r^2}, \] showing that we have exponential convergence.

\sect{Hermite-Pad\'e approximation}
Hermite-Pad\'e approximation is simultaneous rational
approximation to a vector of $r$ functions $f_1,f_2,\ldots,f_r$,
which are all given as Taylor series around a point $a \in
\mathbb{C}$ and for which we require interpolation conditions at
$a$. We will restrict our attention to Hermite-Pad\'e
approximation around infinity and impose interpolation conditions
at infinity.
\subsection{Definition}
Suppose we are given $r$ functions with Laurent expansions
\[     f_j(z) = \sum_{k=0}^\infty \frac{c_{k,j}}{z^{k+1}}, \qquad
j=1,2,\ldots,r. \]
 There are basically two different types of
Hermite-Pad\'e approximation. First we will need multi-indices
$\vec{n}=(n_1,n_2,\ldots,n_r) \in \mathbb{N}^r$ and their size
$|\vec{n}| = n_1+n_2+\cdots+n_r$.

\begin{definition}[Type I]
Type I Hermite-Pad\'e approximation to the vector
$(f_1,\ldots,f_r)$ near infinity consists of finding a vector  
$(A_{\vec{n},1},\ldots,A_{\vec{n},r})$ of polynomials and a
polynomial $B_{\vec{n}}$, with $A_{\vec{n},j}$ of degree $\leq
n_j-1$, such that
\begin{equation}  \label{eq:typeI}
     \sum_{j=1}^r A_{\vec{n},j}(z) f_j(z) - B_{\vec{n}}(z) =
     \O\left(\frac{1}{z^{|\vec{n}|}} \right), \qquad z\to \infty.
\end{equation}
\end{definition}

In type I Hermite-Pad\'e approximation one wants to approximate a
linear combination (with polynomial coefficients) of the $r$
functions by a polynomial. This is often done for the vector of
functions $f,f^2,\ldots,f^r$, where $f$ is a given function. The
solution of the equation
\[   \sum_{j=1}^r A_{\vec{n},j}(z) \hat{f}^j(z) - B_{\vec{n}}(z) = 0 \]
is an algebraic function which gives an algebraic approximant
$\hat{f}$ for the function $f$.

\begin{definition}[Type II]
Type II Hermite-Pad\'e approximation to the vector
$(f_1,\ldots,f_r)$ near infinity consists of finding a polynomial
$P_{\vec{n}}$ of degree $\leq |\vec{n}|$ and polynomials
$Q_{\vec{n},j}$ $(j=1,2,\ldots,r)$ such that
\begin{eqnarray} \label{eq:typeII}
  P_{\vec{n}}(z) f_1(z) - Q_{\vec{n},1}(z) & = & \O \left(
  \frac{1}{z^{n_1+1}} \right), \qquad z \to \infty \nonumber \\
    & \vdots &  \\
P_{\vec{n}}(z) f_r(z) - Q_{\vec{n},r}(z) & = & \O \left(
  \frac{1}{z^{n_r+1}} \right), \qquad z \to \infty. \nonumber
\end{eqnarray}
\end{definition}

Type II Hermite-Pad\'e approximation therefore corresponds to an
approximation of each function $f_j$ separately by rational
functions \textit{with a common denominator} $P_{\vec{n}}$.
Combinations of type I and type II Hermite-Pad\'e approximation
are also possible.

\subsection{Orthogonality}
When we consider $r$ Markov functions
\[   f_j(z) = \int_{a_j}^{b_j} \frac{d\mu_j(x)}{z-x}, \qquad
j=1,2,\ldots,r, \] then Hermite-Pad\'e approximation corresponds
again to certain orthogonality conditions.

First consider type I approximation. Multiply (\ref{eq:typeI}) by
$z^k$ and integrate over a contour $\Gamma$ encircling all the
intervals $[a_j,b_j]$ in the positive direction. Then
\begin{multline*} 
  \frac{1}{2\pi i} \int_\Gamma \left( \sum_{j=1}^r z^k
A_{\vec{n},j}(z)f_j(z)\right) \, dz - \frac{1}{2\pi i} \int_\Gamma z^k
B_{\vec{n}}(z)\, dz \\
= \sum_{\ell=|\vec{n}|}^\infty b_{\vec{n},\ell}
\frac{1}{2\pi i} \int_\Gamma z^{k-\ell}\, dz , 
\end{multline*} 
where the $b_{\vec{n},\ell}$ are the coefficients of the Laurent expansion
of the left hand side in (\ref{eq:typeI}).
Cauchy's theorem implies
\[  \frac{1}{2\pi i} \int_\Gamma z^k
B_{\vec{n}}(z)\, dz = 0. \] Furthermore, there is only a
contribution on the right hand side when $\ell = k+1$, so when $k
\leq |\vec{n}|-2$, then none of the terms in the infinite sum has
a contribution. Therefore we see that
\[  \frac{1}{2\pi i} \int_\Gamma \left( \sum_{j=1}^r z^k
A_{\vec{n},j}(z)f_j(z)\right) \, dz = 0, \qquad 0 \leq k \leq
|\vec{n}|-2.
\] 
Now each $f_j$ is a Markov function, so by changing the order
of integration we get
\[  \frac{1}{2\pi i} \int_\Gamma  z^k
A_{\vec{n},j}(z)f_j(z) \, dz =  \int_{a_j}^{b_j} d\mu_j(x) \
\frac{1}{2\pi i} \int_\Gamma \frac{z^k A_{\vec{n},j}(z)}{z-x}\,
dz. \] Since $\Gamma$ is a contour encircling $[a_j,b_j]$ we have
that
\[  \frac{1}{2\pi i} \int_\Gamma \frac{z^k A_{\vec{n},j}(z)}{z-x}\,
dz = x^k A_{\vec{n},j}(x), \] so that we get the following
orthogonality conditions
\begin{equation}  \label{eq:typeIortho}
   \sum_{j=1}^r \int_{a_j}^{b_j} x^k A_{\vec{n},j}(x)\, d\mu_j(x)
   = 0, \qquad k=0,1,\ldots,|\vec{n}|-2.
\end{equation}
These are $|\vec{n}|-1$ linear and homogeneous equations for the
$|\vec{n}|$ coefficients of the $r$ polynomials $A_{\vec{n},j}$
$(j=1,2,\ldots,r)$, so that we can determine these polynomials up
to a multiplicative factor, provided that the rank of the matrix
in this system is $|\vec{n}|-1$. If the solution is unique (up to
a multiplicative factor), then we say that $\vec{n}$ is a \dword{normal
index for type I}. One can show that this is equivalent to the
condition that the degree of each $A_{\vec{n},j}$ is exactly
$n_j-1$. Once the polynomial vector
$(A_{\vec{n},1},\ldots,A_{\vec{n},r})$ is determined, we can also
find the remaining polynomial $B_{\vec{n}}$ which is given by
\begin{equation}  \label{eq:B}
  B_{\vec{n}}(z) = \sum_{j=1}^r \int_{a_j}^{b_j}
  \frac{A_{\vec{n},j}(z)-A_{\vec{n},j}(x)}{z-x}\, d\mu_j(x).
\end{equation}
Indeed, with this definition of $B_{\vec{n}}$ we have
\begin{equation}  \label{eq:typeIerror}
 \sum_{j=1}^r A_{\vec{n},j}(z) f_j(z) -
B_{\vec{n}}(z) = \sum_{j=1}^r \int_{a_j}^{b_j}
\frac{A_{\vec{n},j}(x)}{z-x}\, d\mu_j(x).
\end{equation}
If we use the expansion
\[  \frac{1}{z-x} = \sum_{k=0}^\infty \frac{x^k}{z^{k+1}}, \]
then the right hand side is
\[   \sum_{k=0}^\infty \frac{1}{z^{k+1}} \sum_{j=1}^r
\int_{a_j}^{b_j} x^k A_{\vec{n},j}(x)\, d\mu_j(x), \] and the
orthogonality conditions (\ref{eq:typeIortho}) show that the sum
over $k$ starts with $k=|\vec{n}|-1$, hence the right hand side is
$\O(z^{-|\vec{n}|})$, which is the order given in the definition
of type I Hermite-Pad\'e approximation.

Next we consider type II approximation. Multiply (\ref{eq:typeII})
by $z^k$ and integrate over a contour $\Gamma$ encircling all the
intervals $[a_j,b_j]$. Then
\begin{multline*}  
\frac{1}{2\pi i} \int_\Gamma z^k P_{\vec{n}}(z) f_j(z)\, dz -
   \frac{1}{2\pi i} \int_\Gamma z^k Q_{\vec{n},j}(z)\, dz \\
   =    \sum_{\ell=n_j+1}^\infty b_{\vec{n},j,\ell} \frac{1}{2\pi i} \int_\Gamma
   z^{k-\ell}\, dz, 
\end{multline*}
where the $b_{\vec{n},j,\ell}$ are the coefficients in the Laurent expansion
of the left hand side of (\ref{eq:typeII}).
 Cauchy's theorem gives
\[  \frac{1}{2\pi i} \int_\Gamma z^k Q_{\vec{n},j}(z)\, dz = 0, \]
and on the right hand side we only have a contribution when
$\ell=k+1$. So for $k \leq n_j-1$ none of the terms in the
infinite sum contribute. Hence
\[  \frac{1}{2\pi i} \int_\Gamma z^k P_{\vec{n}}(z) f_j(z)\, dz =
0, \qquad 0 \leq k \leq n_j-1. \] Interchanging the order of
integration on the left hand side gives the orthogonality
conditions
\begin{eqnarray}  \label{eq:typeIIortho}
    \int_{a_1}^{b_1} x^k P_{\vec{n}}(x) \, d\mu_1(x) & = & 0,
    \qquad k=0,1,\ldots,n_1-1, \nonumber \\
     & \vdots & \\
    \int_{a_r}^{b_r} x^k P_{\vec{n}}(x) \, d\mu_r(x) & = & 0,
    \qquad k=0,1,\ldots,n_r-1. \nonumber
\end{eqnarray}
This gives $|\vec{n}|$ linear and homogeneous equations for the
$|\vec{n}|+1$ coefficients of $P_{\vec{n}}$, hence we can obtain
the polynomial $P_{\vec{n}}$ up to a multiplicative factor,
provided the matrix of coefficients has rank $|\vec{n}|$. In that
case we call the index $\vec{n}$ \dword{normal for type II}. One can show
that this is equivalent to the condition that the degree of
$P_{\vec{n}}$ be exactly $|\vec{n}|$. Once the polynomial
$P_{\vec{n}}$ is determined, we can obtain the polynomials
$Q_{\vec{n},j}$ by
\begin{equation} \label{eq:Qmul}
   Q_{\vec{n},j}(z) = \int_{a_j}^{b_j}
   \frac{P_{\vec{n}}(z)-P_{\vec{n}}(x)}{z-x}\,d\mu_j(x).
\end{equation}
Indeed, with this expression for $Q_{\vec{n},j}$ we have
\begin{equation} \label{eq:typeIIerror}
   P_{\vec{n}}(z)f_j(z) - Q_{\vec{n},j}(z) = \int_{a_j}^{b_j}
   \frac{P_{\vec{n}}(x)}{z-x}\, d\mu_j(x),
\end{equation}
and if we expand $1/(z-x)$, then the right hand side is of the
form
\[  \sum_{k=0}^\infty \frac{1}{z^{k+1}} \int_{a_j}^{b_j} x^k
P_{\vec{n}}(x) \, d\mu_j(x), \] and the orthogonality conditions
(\ref{eq:typeIIortho}) show that the infinite sum starts at
$k=n_j$, which gives an expression of $\O(z^{-n_j-1})$, which is
exactly what is required for type II Hermite-Pad\'e approximation.

\subsection{Angelesco systems}
Angelesco \cite{Ang} introduced an interesting system about which more can be said.
\begin{definition}
An Angelesco system $(f_1,f_2,\ldots,f_r)$ consists of $r$ Markov
functions for which the intervals $(a_j,b_j)$ are pairwise
disjoint.
\end{definition}

All multi-indices are normal for type II in an Angelesco system.
We will prove this by showing that the multiple orthogonal
polynomial $P_{\vec{n}}$ has degree exactly equal to $|\vec{n}|$.
In fact more is true, namely:

\begin{theorem}        \label{thm:ang}
If $(f_1,\ldots,f_r)$ is an Angelesco system with measures $\mu_j$
that have infinitely many points in their support, then
$P_{\vec{n}}$ has $n_j$ simple zeros on $(a_j,b_j)$ for
$j=1,\ldots,r$.
\end{theorem}

\begin{proof}
Let $x_1,\ldots,x_m$ be the sign changes of $P_{\vec{n}}$ on
$(a_j,b_j)$. Suppose that $m < n_j$ and let
$\pi_m(x):=(x-x_1)\cdots(x-x_m)$. Then $P_{\vec{n}}\pi_m$ does not
change sign on $[a_j,b_j]$. Since the support of $\mu_j$ has
infinitely many points, we have
\[  \int_{a_j}^{b_j} P_{\vec{n}}(x)\pi_m(x)\, d\mu_j(x) \neq 0. \]
However, the orthogonality (\ref{eq:typeIIortho}) implies that
$P_{\vec{n}}$ is orthogonal to all polynomials of degree $\leq
n_j-1$ with respect to the measure $\mu_j$ on $[a_j,b_j]$, so that
the integral is zero. This contradiction implies that $m \geq
n_j$, and hence $P_{\vec{n}}$ has at least $n_j$ zeros on
$(a_j,b_j)$. This holds for every $j$, and since the intervals
$(a_j,b_j)$ are disjoint this gives at least $|\vec{n}|$ zeros on
the real line. But the degree of $P_{\vec{n}}$ is $\leq
|\vec{n}|$, hence $P_{\vec{n}}$ has exactly $n_j$ simple zeros on
$(a_j,b_j)$.
\end{proof}

The polynomial $P_{\vec{n}}$ can therefore be factored as
\[   P_{\vec{n}}(x) = q_{n_1}(x)q_{n_2}(x) \cdots q_{n_r}(x), \]
where each $q_{n_j}$ is a polynomial of degree $n_j$ with its
zeros on $(a_j,b_j)$. The orthogonality (\ref{eq:typeIIortho})
then gives \begin{equation}  \label{eq:q}
   \int_{a_j}^{b_j} x^k q_{n_j}(x)\
\prod_{i\neq j} q_{n_i}(x)\, d\mu_j(x) = 0, \qquad k=0,1,\ldots,
n_j-1.
\end{equation}
The product $\prod_{i\neq j} q_{n_i}(x)$ does not change sign on
$(a_j,b_j)$, hence (\ref{eq:q}) shows that $q_{n_j}$ is an
ordinary orthogonal polynomial of degree $n_j$ on the interval
$[a_j,b_j]$ with respect to the measure $\prod_{i\neq j}
|q_{n_i}(x)|\, d\mu_j(x)$. The measure  depends on the multi-index
$\vec{n}$.

\subsection{Algebraic Chebyshev systems}
A Chebyshev system $\{ \varphi_1,\ldots,\varphi_n \}$ on $[a,b]$
is a linearly independent system of
$n$ functions such that 
every nontrivial
linear combination $\sum_{k=1}^n a_k\varphi_k$ has at most $n-1$
zeros on $[a,b]$. This is equivalent to the condition that
\[  \det \begin{pmatrix}
     \varphi_1(x_1) & \varphi_1(x_2) & \cdots & \varphi_1(x_n) \\
     \varphi_2(x_1) & \varphi_2(x_2) & \cdots & \varphi_2(x_n) \\
     \vdots & \vdots & \cdots & \vdots \\
     \varphi_n(x_1) & \varphi_n(x_2) & \cdots & \varphi_n(x_n)
     \end{pmatrix} \neq 0
\]
for every choice of $n$ distinct points $x_1,\ldots,x_n \in
[a,b]$. Indeed, when $x_1,\ldots,x_n$ are such that the
determinant is zero, then there is a linear combination of the
rows that gives a zero row, but this means that for this linear
combination $\sum_{k=1}^n a_k \varphi_k$ has zeros at
$x_1,\ldots,x_n$, giving $n$ zeros, which is not allowed.

\begin{definition}
A system $(f_1,\ldots,f_r)$ is an algebraic Chebyshev system (AT
system) for the index $\vec{n}$ if each $f_j$ is a Markov function
on the same interval $[a,b]$ with a measure $w_j(x)\,d\mu(x)$,
where $\mu$ has infinite support and the $w_j$ are such that
\begin{multline} \label{eq:Chebsys}
 \{w_1,xw_1,\ldots,x^{n_1-1}w_1, w_2, xw_2,\ldots,x^{n_2-1}w_2,
\ldots, \\
w_r,xw_r,\ldots,x^{n_r-1}w_r\}
\end{multline}
is a Chebyshev system on $[a,b]$.
\end{definition}

\begin{theorem}
Suppose $\vec{n}$ is a multi-index such that $(f_1,\ldots,f_r)$ is
an AT system on $[a,b]$ for every index $\vec{m}$ for which $m_j
\leq n_j$ $(1 \leq j \leq r)$. Then $P_{\vec{n}}$ has $|\vec{n}|$
zeros on $(a,b)$ and hence $\vec{n}$ is a normal index for type
II.
\end{theorem}

\begin{proof}
Let $x_1,\ldots,x_m$ be the sign changes of $P_{\vec{n}}$ on
$(a,b)$ and suppose that $m < |\vec{n}|$. We can then find a
multi-index $\vec{m}$ such that $|\vec{m}|=m$ and $m_j \leq n_j$
for every $1 \leq j \leq r$ and $m_k < n_k$ for some $1 \leq k
\leq r$. Consider the interpolation problem where we want to find
a function
\[   L(x) = \sum_{j=1}^r q_j(x)w_j(x), \]
where $q_j$ is a polynomial of degree $m_j-1$ if $j\neq k$ and
$q_k$ a polynomial of degree $m_k$ that satisfies
\begin{eqnarray*}
 L(x_j) & = & 0, \qquad j=1,...,m, \\
 L(x_0) & = & 1, \qquad \textrm{for some other point $x_0 \in
 [a,b]$}.
\end{eqnarray*}
The function $L$ is a linear combination of
\[  \{w_1,xw_1,\ldots,x^{m_1-1}w_1,\ldots, w_k, xw_k,\ldots,x^{n_k}w_k, \ldots
  w_r,xw_r,\ldots,x^{m_r-1}w_r\} \]
and this is, by assumption, a Chebyshev system.
This interpolation problem has a unique solution since it
involves a Chebyshev system of basis functions. The function $L$
has, by construction, $m$ zeros and the Chebyshev system has $m+1$
basis functions, so $L$ can have at most $m$ zeros on $[a,b]$ and
each zero is a sign change (see, e.g., \cite[pp.~20--21]{KS}).
Hence $P_{\vec{n}}L$ does not change
sign on $[a,b]$.
 Since $\mu$ has infinite support, we
thus have
\[  \int_a^b L(x)P_{\vec{n}}(x)\, d\mu(x) \neq 0. \]
But the orthogonality (\ref{eq:typeIIortho}) gives
\[  \int_a^b q_{j}(x)P_{\vec{n}}(x) w_j(x)\, d\mu(x) = 0, \qquad
j=1,2,\ldots,r, \] and this contradiction implies that
$P_{\vec{n}}$ has $|\vec{n}|$ simple zeros on $(a,b)$.
\end{proof}

We have a similar result for type I Hermite-Pad\'e approximation:

\begin{theorem}
Suppose $\vec{n}$ is a multi-index such that $(f_1,\ldots,f_r)$ is
an AT system on $[a,b]$ for every index $\vec{m}$ for which $m_j
\leq n_j$ $(1 \leq j \leq r)$. Then $\sum_{j=1}^r
A_{\vec{n},j}w_j$ has $|\vec{n}|-1$ zeros on $(a,b)$ and
$\vec{n}$ is a normal index for type I.
\end{theorem}

\begin{proof}
Let $x_1,\ldots,x_m$ be the sign changes of $\sum_{j=1}^r
A_{\vec{n},j}w_j$ on $(a,b)$ and suppose that $m < |\vec{n}|-1$.
Let $\pi_m$ be the monic polynomial with these points as zeros. Then
$\pi_m \sum_{j=1}^r A_{\vec{n},j}w_j$ does not change sign on
$[a,b]$ and hence
\[  \int_a^b \pi_m(x) \sum_{j=1}^r
A_{\vec{n},j}(x) w_j(x)\, d\mu(x) \neq 0. \] But the orthogonality
conditions (\ref{eq:typeIortho}) indicate that this integral is
zero. This contradiction implies that $m \geq |\vec{n}|-1$. The
sum $\sum_{j=1}^r A_{\vec{n},j}w_j$ is a linear combination of the
Chebyshev system (\ref{eq:Chebsys}), hence it has at most
$|\vec{n}|-1$ zeros on $[a,b]$. Therefore we see that $m =
|\vec{n}|-1$. To see that the index $\vec{n}$ is normal for type
I, we assume that for some $k$ with $1 \leq k \leq r$ the degree
of $A_{\vec{n},k}$ is less than $n_k-1$. Then $\sum_{j=1}^r
A_{\vec{n},j}w_j$ is a linear combination of the Chebyshev
system (\ref{eq:Chebsys}) from which the function $x^{n_k-1}w_k$
is removed. This is still a Chebyshev system by assumption, and
hence this linear combination has at most $|\vec{n}|-2$ zeros on
$[a,b]$. But this contradicts our previous observation that it has
$|\vec{n}|-1$ zeros. Therefore every $A_{\vec{n},j}$ has degree
exactly $n_j-1$, so that the index $\vec{n}$ is normal.
\end{proof}

\subsection{Nikishin systems}
A special construction, suggested by Nikishin \cite{Nik}, gives an AT system
that can be handled in some detail. The construction is by
induction. A \dword{Nikishin system of order 1} is a Markov function
$f_{1,1}$ for a measure $\mu_1$ on the interval $[a_1,b_1]$. A
\dword{Nikishin system of order 2} is a vector of Markov functions
$(f_{1,2}, f_{2,2})$ on $[a_2,b_2]$ such that
\[    f_{1,2}(z) = \int_{a_2}^{b_2} \frac{d\mu_2(x)}{z-x}, \quad
      f_{2,2}(z) = \int_{a_2}^{b_2} f_{1,1}(x)
      \frac{d\mu_2(x)}{z-x}, \]
where $f_{1,1}$ is a Nikishin system of order 1 on $[a_1,b_1]$ and
$(a_1,b_1) \cap (a_2,b_2) = \emptyset$. In general we have
\begin{definition}
A \dword{Nikishin system of order $r$} consists of $r$ Markov functions
$(f_{1,r},\ldots,f_{r,r})$ on $[a_r,b_r]$ such that
 \begin{eqnarray}
f_{1,r}(z) & = &  \int_{a_r}^{b_r} \frac{d\mu_r(x)}{z-x}, \\
f_{j,r}(z) & = & \int_{a_r}^{b_r} f_{j-1,r-1}(x)
\frac{d\mu_r(x)}{z-x}, \qquad j=2,\ldots,r,
\end{eqnarray}
where $(f_{1,r-1},\ldots,f_{r-1,r-1})$ is a Nikishin system of
order $r-1$ on $[a_{r-1},b_{r-1}]$ and $(a_r,b_r) \cap
(a_{r-1},b_{r-1}) = \emptyset$.
\end{definition}

For a Nikishin system of order $r$ one knows that the
multi-indices $\vec{n}$ with $n_1 \geq n_2 \geq \cdots \geq n_r$
are normal (the system is an AT-system for these indices), but it
is an open problem whether every multi-index is normal (for $r>2$;
for $r=2$ it has been proved that every multi-index is normal).

What can be said about type II Hermite-Pad\'e approximation for
$r=2$? Recall (\ref{eq:typeIIerror}) for the function $f_{1,2}$:
\[   P_{n_1,n_2}(y) f_{1,2}(y) - Q_{n_1,n_2;1}(y) = \int_{a_2}^{b_2}
   \frac{P_{n_1,n_2}(x)}{y-x}\,d\mu_2(x). \]
Multiply both sides by $y^k$, with $k\leq n_1$. Then the right
hand side is
\begin{multline*}  
  \int_{a_2}^{b_2}
   \frac{y^kP_{n_1,n_2}(x)}{y-x}\,d\mu_2(x) \\
 = \int_{a_2}^{b_2}
   \frac{(y^k-x^k)P_{n_1,n_2}(x)}{y-x}\,d\mu_2(x) + \int_{a_2}^{b_2}
   \frac{x^kP_{n_1,n_2}(x)}{y-x}\,d\mu_2(x). 
\end{multline*}
Clearly $(y^k-x^k)/(y-x)$ is a polynomial in $x$ of degree $k-1
\leq n_1-1$ hence the first integral on the right vanishes because
of the orthogonality (\ref{eq:typeIIortho}). Integrate over the
variable $y \in [a_1,b_1]$ with respect to the measure $\mu_1$. Then
we find for $k \leq n_1$
\begin{multline*}
  \int_{a_1}^{b_1} [P_{n_1,n_2}(y) f_{1,2}(y) - Q_{n_1,n_2;1}(y)]
y^k \, d\mu_1(y) \\
= \int_{a_1}^{b_1} \int_{a_2}^{b_2} \frac{x^k
P_{n_1,n_2}(x)}{y-x} \, d\mu_2(x)\,d\mu_1(y). 
\end{multline*}
Change the order
of integration on the right hand side. Then
\begin{multline*}
  \int_{a_1}^{b_1} [P_{n_1,n_2}(y) f_{1,2}(y) - Q_{n_1,n_2;1}(y)]
y^k \, d\mu_1(y) \\
= - \int_{a_2}^{b_2}  x^k P_{n_1,n_2}(x)
f_{1,1}(x) \, d\mu_2(x) 
\end{multline*} 
and this is zero for $k \leq n_2-1$.
Hence if $n_2 \leq n_1+1$ then the expression $P_{n_1,n_2}
f_{1,2} - Q_{n_1,n_2;1}$ is orthogonal to all polynomials of
degree $\leq n_2-1$ on $[a_1,b_1]$. This implies that
$P_{n_1,n_2} f_{1,2} - Q_{n_1,n_2;1}$ has at least $n_2$
zeros on $(a_1,b_1)$ using an argument similar to what we have been
using earlier. Let $R_{n_2}$ be the monic polynomial with $n_2$ of
these zeros on $(a_1,b_1)$. Then $[P_{n_1,n_2} f_{1,2} -
Q_{n_1,n_2;1}]/R_{n_2}$ is an analytic function on
$\mathbb{C} \setminus [a_2,b_2]$, which has the representation
\[  \frac{P_{n_1,n_2}(y) f_{1,2}(y) - Q_{n_1,n_2;1}(y)}{R_{n_2}(y)}
   = \frac{1}{R_{n_2}(y)} \int_{a_2}^{b_2}
   \frac{P_{n_1,n_2}(x)}{y-x}\, d\mu_2(x). \]
Multiply both sides by $y^k$ and integrate over a contour $\Gamma$
encircling the interval $[a_2,b_2]$ in the positive direction, but
with all the zeros of $R_{n_2}$ outside $\Gamma$. Then
\begin{multline*}
 \frac{1}{2\pi i} \int_\Gamma y^k \frac{P_{n_1,n_2}(y) f_{1,2}(y) -
Q_{n_1,n_2;1}(y)}{R_{n_2}(y)}\, dy  \\
   = \frac{1}{2\pi i} \int_\Gamma \frac{y^k}{R_{n_2}(y)}
   \frac{P_{n_1,n_2}(x)}{y-x}\, d\mu_2(x)\, dy. 
\end{multline*}
If we interchange the order of integration on the right hand side
and use Cauchy's theorem, then this gives the integral
\[ \int_{a_2}^{b_2} x^k P_{n_1,n_2}(x)\,
\frac{d\mu_2(x)}{R_{n_2}(x)}. \] 
By the interpolation condition
(\ref{eq:typeII}), the integrand on the left hand side is of the
order $\O(y^{k-n_1-n_2-1})$, so if we use Cauchy's theorem for the
exterior of $\Gamma$, then we see that the integral vanishes for $k \leq
n_1+n_2-1$. Hence we get
\begin{equation} \label{eq:typeIInikortho}
\int_{a_2}^{b_2} x^k P_{n_1,n_2}(x)\, \frac{d\mu_2(x)}{R_{n_2}(x)}
= 0, \qquad k=0,1,\ldots, n_1+n_2-1.
\end{equation}
This shows that $P_{n_1,n_2}$ is an ordinary orthogonal polynomial
on $[a_2,b_2]$ with respect to the measure $d\mu_2/R_{n_2}$.
Observe that $(a_1,b_1) \cap (a_2,b_2) = \emptyset$ implies that
$R_{n_2}$ does not change sign on $[a_2,b_2]$. Finally we have
\begin{eqnarray*} \int_{a_2}^{b_2} \frac{P_{n_1,n_2}^2(x)}{y-x} \,
\frac{d\mu_2(x)}{R_{n_2}(x)} & = & \int_{a_2}^{b_2}
P_{n_1,n_2}(x)\frac{P_{n_1,n_2}(x)-P_{n_1,n_2}(y)}{y-x} \,
\frac{d\mu_2(x)}{R_{n_2}(x)}  \\
 & & +\ P_{n_1,n_2}(y) \int_{a_2}^{b_2}
\frac{P_{n_1,n_2}(x)}{y-x} \, \frac{d\mu_2(x)}{R_{n_2}(x)} \\
 & = & P_{n_1,n_2}(y) \int_{a_2}^{b_2}
\frac{P_{n_1,n_2}(x)}{y-x} \, \frac{d\mu_2(x)}{R_{n_2}(x)},
\end{eqnarray*}
since $[P_{n_1,n_2}(y)-P_{n_1,n_2}(x)]/(y-x)$ is a polynomial in
$x$ of degree $n_1+n_2-1$ and because of the orthogonality
(\ref{eq:typeIInikortho}). Hence
\begin{equation} \label{eq:typeIInikerror}
P_{n_1,n_2}(y) f_{1,2}(y) - Q_{n_1,n_2;1}(y) =
\frac{R_{n_2}(y)}{P_{n_1,n_2}(y)} \int_{a_2}^{b_2}
\frac{P_{n_1,n_2}^2(x)}{y-x} \, \frac{d\mu_2(x)}{R_{n_2}(x)} .
\end{equation}
Both sides of the equation have zeros at the zeros of $R_{n_2}$,
but there will not be any other zeros on $[a_1,b_1]$ since the
integral on the right hand side has constant sign.

\subsection{Asymptotic properties and convergence}
We restrict ourselves to the case $r=2$, but the general case
$r>1$ can be treated in a similar way (with a bit more work). The
asymptotic properties of the multiple orthogonal polynomials and
the convergence of the Hermite-Pad\'e approximants are handled by
trying to put everything into terms of ordinary orthogonal
polynomials.
\subsubsection{Angelesco systems}
The type II multiple orthogonal polynomial can be factored as
$P_{n_1,n_2} = q_{n-1}q_{n-2}$,
where $q_{n_1}$ has $n_1$
zeros on $(a_1,b_1)$ and $q_{n_2}$ has $n_2$ zeros on $(a_2,b_2)$.
 From (\ref{eq:typeIIerror}) we get
\[  f_1(z) - \frac{Q_{n_1,n_2;1}(z)}{P_{n_1,n_2}(z)} =
\frac{1}{q_{n_1}(z)q_{n_2}(z)} \int_{a_1}^{b_1}
\frac{q_{n_1}(x)}{z-x} \, q_{n_2}(x)\, d\mu_1(x). \] We saw that
$q_{n_1}$ is an orthogonal polynomial of degree $n_1$ on
$[a_1,b_1]$ for the measure $|q_{n_2}(x)|\, d\mu_1(x)$, so we can
write
\[ \int_{a_1}^{a_2}
\frac{q_{n_1}(x)}{z-x} \, q_{n_2}(x)\, d\mu_1(x) =
\frac{1}{q_{n_1}(z)} \int_{a_1}^{b_1} \frac{q_{n_1}^2(x)}{z-x} \,
q_{n_2}(x)\, d\mu_1(x) \] as we did earlier in Section 1.6. This
gives
\[   f_1(z) - \frac{Q_{n_1,n_2;1}(z)}{P_{n_1,n_2}(z)} =
\frac{1}{q_{n_1}^2(z)q_{n_2}(z)} \int_{a_1}^{a_2}
\frac{q_{n_1}^2(x)}{z-x} \, q_{n_2}(x)\, d\mu_1(x). \] 
 From here
we get the estimate
\[ \left| f_1(z) - \frac{Q_{n_1,n_2;1}(z)}{P_{n_1,n_2}(z)} \right| \leq
\frac{1}{|q_{n_1}(z)|^2 |q_{n_2}(z)|} \frac{1}{d_1}
\int_{a_1}^{b_1} q_{n_1}^2(x) \, |q_{n_2}(x)|\, d\mu_1(x),
\]
where $d_1$ is the distance between $z$ and $[a_1,b_1]$. If
$P_{n_1,n_2}$ is normalized so that it is monic, then we can take
both $q_{n_1}$ and $q_{n_2}$ monic and we get
\[  \left| f_1(z) - \frac{Q_{n_1,n_2;1}(z)}{P_{n_1,n_2}(z)} \right| \leq
\frac{1}{d_1 \gamma_{n_1;1}^2|q_{n_1}(z)|^2 |q_{n_2}(z)|} ,
\]
where
\begin{eqnarray}  
  \frac{1}{\gamma_{n_1;1}^2} &=& \int_{a_1}^{b_1}
q_{n_1}^2(x) \, |q_{n_2}(x)|\, d\mu_1(x) \nonumber \\
  &=& \min_{\pi_{n_1}(x)=x^{n_1}+\cdots} \int_{a_1}^{b_1} \pi_{n_1}^2(x) \, |q_{n_2}(x)|\,
  d\mu_1(x).         \label{eq:gamma1}
\end{eqnarray}
A similar reasoning holds for the rational approximation to $f_2$
and gives
\[  \left| f_2(z) - \frac{Q_{n_1,n_2;2}(z)}{P_{n_1,n_2}(z)} \right| \leq
\frac{1}{d_2 \gamma_{n_2;2}^2|q_{n_2}(z)|^2 |q_{n_1}(z)|} ,
\]
where $d_2$ is the distance of $z$ to $[a_2,b_2]$ and
\begin{eqnarray} 
 \frac{1}{\gamma_{n_1;2}^2} &=& \int_{a_2}^{b_2}
q_{n_2}^2(x) \, |q_{n_1}(x)|\, d\mu_2(x) \nonumber \\
 & =& \min_{\pi_{n_2}(x)=x^{n_2}+\cdots} \int_{a_2}^{b_2} \pi_{n_2}^2(x) \, |q_{n_1}(x)|\,
  d\mu_2(x).   \label{eq:gamma2}
\end{eqnarray}
 The convergence of these rational approximants is therefore
given in terms of the asymptotic behavior of $|q_{n_1}(z)|$,
$|q_{n_2}(z)|$ and the constants $\gamma_{n_1;1}$ and
$\gamma_{n_2;2}$. These polynomials (and their zeros) interact
with each other: the polynomial $q_{n_1}$ is an orthogonal
polynomial for a measure that contains $q_{n_2}$ as a factor, and
$q_{n_2}$ is an orthogonal polynomial for a measure that contains
$q_{n_1}$ as a factor. Let
\[   \nu_{n_1;1} := \frac{1}{n_1} \sum_{j=1}^{n_1}
\delta_{x_{j,n_1}}, \quad \nu_{n_2;2} := \frac{1}{n_2}
\sum_{j=1}^{n_2} \delta_{y_{j,n_2}}, \] 
where $x_{j,n_1}$ are the
zeros of $q_{n_1}$ and $y_{j,n_2}$ are the zeros of $q_{n_2}$.
Then
$(\nu_{n_1;1})$ is a sequence of probability measures on
$[a_1,b_1]$ and $(\nu_{n_2;2})$ is a sequence of probability
measures on $[a_2,b_2]$. Helly's selection principle guarantees
that there are weakly converging subsequences with limits $\nu_1$
on $[a_1,b_1]$ and $\nu_2$ on $[a_2,b_2]$.
 The minimization problems (\ref{eq:gamma1}) and
(\ref{eq:gamma2}) lead to an extremal problem in potential theory
for two probability measures. The integral in (\ref{eq:gamma1}) is
approximately of the form
\[  \int_{a_1}^{b_1} \exp \left[ -2n_1 U(x;\nu_1) -n_2 U(x;\nu_2)
\right]\, d\mu_1(x) \] where $U(x;\nu)$ is the logarithmic potential of $\nu$
\[  U(x;\nu) = \int \log \frac{1}{|x-y|} \, d\nu(y), \]
and the integral in (\ref{eq:gamma2}) is approximately of the form
\[  \int_{a_2}^{b_2} \exp \left[ -2n_2 U(x;\nu_2) -n_1 U(x;\nu_1)
\right]\, d\mu_2(x). \] We want to minimize both integrals over
all pairs of probability measures $(\nu_1,\nu_2)$, where the first
measure is supported on $[a_1,b_1]$ and the second measure on
$[a_2,b_2]$. If $n_1/(n_1+n_2) \to p$ and $n_2/(n_1+n_2) \to q$
(so that $p+q=1$), and if the measures $\mu_1$ and $\mu_2$ are
sufficiently regular (e.g., $\mu_1' > 0$ almost everywhere on
$[a_1,b_1]$ and $\mu_2' >0$ almost everywhere on $[a_2,b_2]$) then
the solution of the extremal problem satisfies
\begin{eqnarray}
    2pU(x;\nu_1) + q U(x;\nu_2) & = &
\ell_1, \qquad x \in
\textrm{supp}(\nu_1) \subset [a_1,b_1], \label{eq:Uang1} \\
    p U(x;\nu_1) + 2q U(x;\nu_2) & = & \ell_2, \qquad x \in
    \textrm{supp}(\nu_2) \subset [a_2,b_2]. \label{eq:Uang2}
\end{eqnarray}
where the $\ell_j$ are constants that act as Lagrange multipliers.
For this extremal problem it is possible that the support of
$\nu_1$ is not the full interval $[a_1,b_1]$ and the support of
$\nu_2$ can be a subset of $[a_2,b_2]$. This is a consequence of
the interaction: the zeros of $q_{n_1}$ are repelling the zeros of
$q_{n_2}$ and vice versa. The variational conditions
(\ref{eq:Uang1})--(\ref{eq:Uang2}) have to be supplemented with
\begin{eqnarray*}
    2pU(x;\nu_1) + q U(x;\nu_2) & \geq &
\ell_1, \qquad x \in
[a_1,b_1] \setminus \textrm{supp}(\nu_1),  \\
    p U(x;\nu_1) + 2q U(x;\nu_2) & \geq & \ell_2, \qquad x \in
     [a_2,b_2] \setminus \textrm{supp}(\nu_2).
\end{eqnarray*}
The Lagrange multipliers $\ell_1,\ell_2$ appear in the asymptotics
of $\gamma_{n_1;1}$ and $\gamma_{n_2;2}$ as
\[  \lim_{n_1+n_2\to \infty} \gamma_{n_1;1}^{2/(n_1+n_2)} = \exp(
\ell_1), \quad \lim_{n_1+n_2\to \infty}
\gamma_{n_2;2}^{2/(n_1+n_2)} = \exp( \ell_2). \]

Our conclusion is that the convergence to first function $f_1$ is determined by
level curves  $C_r= \{ z : \exp[2pU(z;\nu_1)+qU(z;\nu_2)-\ell_1]=r
\}$ with $r < 1$ on which we have \[ \lim_{n_1+n_2\to \infty}
\left| f_1(z) - \frac{Q_{n_1,n_2;1}(z)}{P_{n_1,n_2}(z)}
\right|^{1/(n_1+n_2)} = r
\]
and the convergence to the second function $f_2$ is determined by level curves $D_r =
\{z : \exp[pU(z;\nu_1)+2qU(z;\nu_2)-\ell_2]=r \}$ with $r < 1$ on
which we have \[ \lim_{n_1+n_2\to \infty} \left| f_2(z) -
\frac{Q_{n_1,n_2;2}(z)}{P_{n_1,n_2}(z)} \right|^{1/(n_1+n_2)} = r.
\]
Observe that $\textrm{supp}(\nu_1) \subset C_1$ and
$\textrm{supp}(\nu_2) \subset D_1$, so we don't expect exponential
convergence on these sets. On the remaining part of $[a_1,b_1]$
(and $[a_2,b_2]$) we get values $r\geq 1$, so we get even worse
behavior there. This is caused by the fact that on these parts of
the intervals there will not be enough zeros of the multiple
orthogonal polynomial to simulate the singularities of the
functions $f_1$ and $f_2$.

\subsubsection{Nikishin systems}
The analysis for Nikishin systems is similar but leads to a
different extremal problem for potentials. We now start from
(\ref{eq:typeIInikerror}) which gives
\begin{equation}  \label{eq:HPtypeII}
 \left|  f_{1,2}(y) - \frac{Q_{n_1,n_2;1}(y)}{P_{n_1,n_2}(y)} \right|
\leq \frac{|R_{n_2}(y)|}{|P_{n_1,n_2}(y)|^2} \frac{1}{d_2}
\int_{a_2}^{b_2} P_{n_1,n_2}^2(x) \,
\frac{d\mu_2(x)}{|R_{n_2}(x)|},
\end{equation} where $d_2$ is the
distance from $y$ to $[a_2,b_2]$. Now we have that $P_{n_1,n_2}$
is a (monic) orthogonal polynomial on $[a_2,b_2]$ for the measure
$d\mu_2/|R_{n_2}|$, so we have
\begin{eqnarray}  
  \frac{1}{\gamma_{n_1,n_2}^2} &=& \int_{a_2}^{b_2} P_{n_1,n_2}^2(x) \,
\frac{d\mu_2(x)}{|R_{n_2}(x)|} \nonumber \\
  &=&
\min_{\pi_{n_1+n_2}(x)=x^{n_1+n_2}+\cdots} \int_{a_2}^{b_2}
\pi_{n_1,n_2}^2(x) \, \frac{d\mu_2(x)}{|R_{n_2}(x)|}. \label{eq:gamma12}
\end{eqnarray}
The polynomial $R_{n_2}$ has its zeros on $[a_1,b_1]$ and in fact
is a monic orthogonal polynomial on $[a_1,b_1]$ for the measure
\[   \frac{P_{n_1,n_2}f_{1,2} -
Q_{n_1,n_2;1}}{R_{n_2}}\, d\mu_1 . \] Indeed, we can
verify that
\begin{multline*}
 \int_{a_1}^{b_1} y^k R_{n_2}(y) \frac{P_{n_1,n_2}(y)f_{1,2}(y) -
Q_{n_1,n_2;1}(y)}{R_{n_2}(y)}\, d\mu_1(y) \\ = \int_{a_1}^{b_1}
y^k [P_{n_1,n_2}(y)f_{1,2}(y) - Q_{n_1,n_2;1}(y)]\, d\mu_1(y) = 0,
\end{multline*}
 for $k \leq n_2-1$, since we have seen that the expression
$P_{n_1,n_2}f_{1,2} - Q_{n_1,n_2;1}$ is orthogonal to all
polynomials of degree less than $n_2$ on $[a_1,b_1]$ for the
measure $\mu_1$. The orthogonality measure for $R_{n_2}$ can also
be written as
\[ \frac{P_{n_1,n_2}(y)f_{1,2}(y) -
Q_{n_1,n_2;1}(y)}{R_{n_2}(y)} = \frac{1}{P_{n_1,n_2}(y)}
\int_{a_2}^{b_2} \frac{P_{n_1,n_2}^2(x)}{y-x} \,
\frac{d\mu_2(x)}{R_{n_2}(x)}. \] In this weight we have
\[  \frac{1}{\gamma_{n_1,n_2}^2 C_1} \leq \int_{a_2}^{b_2} \frac{P_{n_1,n_2}^2(x)}{|y-x|} \,
\frac{d\mu_2(x)}{|R_{n_2}(x)|} \leq \frac{1}{\gamma_{n_1,n_2}^2 C_2},
\] where $C_1$ and $C_2$ are the maximum and minimum,
respectively, over the set
\[      \{ |x-y| : x \in [a_2,b_2],\ y \in [a_1,b_1] \}.      \]
So, up to the constants $C_1,C_2$, we have the extremal problem
\begin{eqnarray}  
  \frac{1}{\gamma_{n_2;2}^2} &=& \int_{a_1}^{b_1} R_{n_2}^2(y) \,
  \frac{d\mu_1(y)}{|P_{n_1,n_2}(y)|} \nonumber \\
 & =&  \min_{\pi_{n_2}(y)=y^{n_2}+
  \cdots} \int_{a_1}^{b_1} \pi_{n_2}^2(y) \,
  \frac{d\mu_1(y)}{|P_{n_1,n_2}(y)|}.   \label{eq:gamma2R}
 \end{eqnarray}

 Define the zero distributions
 \[  \nu_{n_1+n_2} := \frac{1}{n_1+n_2} \sum_{j=1}^{n_1+n_2}
 \delta_{x_{j,n_1+n_2}}, \quad \nu_{n_2;2} := \frac{1}{n_2}
 \sum_{j=1}^{n_2} \delta_{y_{j,n_2}}, \]
 where $x_{j,n_1+n_2}$ are the zeros of $P_{n_1,n_2}$ and
 $y_{j,n_2}$ are the zeros of $R_{n_2}$. Then $(\nu_{n_1+n_2})$ is a
 sequence of probability measures on $[a_2,b_2]$ and
 $(\nu_{n_2;2})$ is a sequence of probability measures on
 $[a_1,b_1]$. Helly's selection principle shows that there are
 weakly convergent subsequences with limits $\nu$ and $\nu_2$
 which are supported on $[a_2,b_2]$ and $[a_1,b_1]$ respectively.
 The extremal problems (\ref{eq:gamma12}) and (\ref{eq:gamma2R})
 then lead to an extremal problem in potential theory. The
 integral in (\ref{eq:gamma12}) is approximately
 \[  \int_{a_2}^{b_2} \exp [-2(n_1+n_2) U(x;\nu) + n_2 U(x;\nu_2)]
 \, d\mu_2(x) \]
 and the integral in (\ref{eq:gamma2R}) is approximately
 \[   \int_{a_1}^{b_1} \exp [ -2n_2 U(x;\nu_2) + (n_1+n_2)
 U(x;\nu)] \, d\mu_1(x). \]
 If $n_2/(n_1+n_2) \to q$ and $\mu_i'>0$ almost everywhere on
 $[a_i,b_i]$ $(i=1,2)$, then this gives the variational conditions
 \begin{eqnarray}
    2 U(x;\nu) - q U(x;\nu_2) & = & \ell_1, \qquad x \in
    \textrm{supp}(\nu) \subset [a_2,b_2], \\
    -U(x;\nu) + 2q U(x;\nu_2) & = & \ell_2, \qquad x \in
    \textrm{supp}(\nu_2) \subset [a_1,b_1],
 \end{eqnarray}
 where $\ell_1$ and $\ell_2$ are Lagrange multipliers for which
 \[ \lim_{n_1+n_2\to \infty} \gamma_{n_1,n_2}^{2/(n_1+n_2)} =
 \exp(\ell_1), \quad
   \lim_{n_1+n_2\to \infty} \gamma_{n_2;2}^{2/(n_1+n_2)} =
 \exp(\ell_2). \]

Looking back to (\ref{eq:HPtypeII}) we thus have
\[  \lim_{n_1+n_2\to \infty}
           \left|  f_{1,2}(y) - \frac{Q_{n_1,n_2;1}(y)}{P_{n_1,n_2}(y)}
           \right|^{1/(n_1+n_2)} = r < 1 \]
on level curves $C_r := \{ z : \exp [2U(z;\nu)-qU(z;\nu_2)-\ell_1] = r\}$.

The convergence to the second function $f_{2,2}$ can also be handled but is left
as an advanced exercise for the reader.
\newpage

\sect{Applications}
\subsection{Gauss and simultaneous Gauss quadrature}
Gauss quadrature is directly related to orthogonal polynomials, and hence
to Pad\'e approximation. Here is an approach based on complex analysis.
Suppose $\mu$ is a positive measure on $[a,b]$ and we denote by $f$ the
Markov function for $\mu$,
\[  f(z) = \int_a^b \frac{d\mu(x)}{z-x}. \]
Let $Q_{n-1}/P_n$ be the Pad\'e approximant to $f$ near infinity.
Then
\[      f(z) - \frac{Q_{n-1}(z)}{P_n(z)} = \O(z^{-2n-1}), \qquad z \to \infty. \]
Multiply both sides by a polynomial $\pi_{2n-1}$ of degree at most $2n-1$,
and integrate along a contour $\Gamma$ encircling the interval $[a,b]$ once in the positive
direction. Then
\[  \frac{1}{2\pi i} \int_\Gamma \pi_{2n-1}(z) f(z) \, dz = \frac{1}{2\pi i} \int_\Gamma
      \pi_{2n-1}(z) \frac{Q_{n-1}(z)}{P_n(z)}\, dz, \]
because the remainder term vanishes after integration, due to Cauchy's theorem for
the outside of $\Gamma$. Interchanging the order of integration on the left hand side and
using the residue theorem on the right hand side shows that for every polynomial
$\pi_{2n-1}$ of degree $\leq 2n-1$ we have
\begin{equation}       \label{eq:Gauss}
     \int_{a}^b \pi_{2n-1}(x)\, d\mu(x) = \sum_{j=1}^n \lambda_{j,n} \pi_{2n-1}(x_{j,n}) ,
\end{equation}
where $\lambda_{j,n}$ is the residue of the Pad\'e approximant at
the zeros $x_{j,n}$ of $P_n$, i.e.,
\[    \lambda_{j,n} = \frac{Q_{n-1}(x_{j,n})}{P_n'(x_{j,n})}. \]
If we take $\pi_{2n-1}(x) = P_n^2(x)/(x-x_{j,n})^2$, then (\ref{eq:Gauss}) gives
\[     \int_a^b   \frac{P_n^2(x)}{(x-x_{j,n})^2} \, d\mu(x) = \lambda_{j,n} [P_n'(x_{j,n})]^2, \]
which shows that $\lambda_{j,n} > 0$ for $j=1,\ldots,n$. These weights $\lambda_{j,n}$ are known as
Christoffel numbers or Gauss quadrature coefficients, the zeros $x_{j,n}$ of $P_n$ are
Gauss quadrature nodes, and (\ref{eq:Gauss}) is the Gauss quadrature formula. Replacing
$\pi_{2n-1}$ by a continuous function $g$ on $[a,b]$, suggests to use the sum
\[            \sum_{j=1}^n \lambda_{j,n} g(x_{j,n}) \]
as an approximation to the integral
\[  \int_a^b g(x)\, d\mu(x). \]
If $[a,b]$ is a finite interval, then every continuous function can be approximated uniformly
by polynomials (Weierstrass), hence the quadrature sum indeed converges to the integral
when the number of nodes $n$ tends to infinity. The positivity of the weights $\lambda_{j,n}$
is needed to get this convergence. The quadrature formula requires $n$ function evaluations
(at the zeros of $P_n$) and is exact for polynomials of degree $\leq 2n-1$, hence on a
linear space of dimension $2n$. The ratio $n/2n=1/2$ is a measure for the efficiency of this
formula.

In a number of applications we need to approximate several integrals of the same function,
but with respect to different measures. The following example comes from \cite{Bo}.
Suppose that $g$ is the spectral distribution of light in the direction of the observer and
$w_1, w_2, w_3$ are weight functions describing the profiles for red, green and blue light. Then
the integrals
\[  \int_0^{2\pi} g(x)w_1(x)\, dx, \quad \int_0^{2\pi} g(x)w_2(x)\, dx, \quad \int_0^{2\pi} g(x)w_3(x)\, dx \]
give the amount of light after passing through the filters for red, green and blue.
In this case we need to approximate three integrals of the same function $g$. We would
like to use as few function evaluations as possible, but the integrals should be
accurate for polynomials $g$ of degree as high as possible. If we use Gauss quadrature
with $n$ nodes for each integral, then we require $3n$ function evaluations and all integrals
will be correct for polynomials of degree $\leq 2n-1$ (a space of dimension
 $2n$). This gives
an efficiency of $3/2$. In fact, with $3n$ function evaluations we can double the dimension
of the space in which the formula is exact. Consider the Markov functions
\[  f_j(z) = \int_0^{2\pi} \frac{w_j(x)\, dx}{z-x}, \qquad j=1,2,3 \]
and the type II Hermite-Pad\'e approximation problem
\[    f_j(z) - \frac{Q_{n,n,n;j}(z)}{P_{n,n,n}(z)} = \O(z^{-4n-1}), \qquad z\to \infty. \]
Now we can multiply by a polynomial $\pi_{4n-1}$ of degree at most $4n-1$, and integrate
along a contour $\Gamma$ encircling $[0,2\pi]$ in the positive direction, to obtain
\begin{equation}  \label{eq:simGauss}
    \int_0^{2\pi} \pi_{4n-1}(x)w_j(x)\, dx = \sum_{k=1}^{3n} \lambda_{k,n;j}
    g(x_{k,n}), \qquad j=1,2,3,
\end{equation}
where $x_{k,n}$ are the zeros of $P_{n,n,n}$ and $\lambda_{k,n;j}$ are the residues
of $Q_{n,n,n;j}/P_{n,n,n}$ at the zero $x_{k,n}$:
\[  \lambda_{k,n;j} =
\frac{Q_{n,n,n;j}(x_{k,n})}{P_{n,n,n}'(x_{k,n})} . \]
Therefore the three integrals will be evaluated exactly by the three sums in
(\ref{eq:simGauss}) for polynomials of degree $\leq 4n-1$. The convergence is
somewhat more difficult to handle, since we do not have a general result that
the quadrature coefficients $\lambda_{k,n;j}$ are positive. The positivity
has to be investigated separately for Angelesco and Nikishin systems. 
See \cite{CV,FIL,FIL2} for finding out more about simultaneous Gauss quadrature.

\subsection{Irrationality and transcendence}
Hermite-Pad\'e approximants were introduced by Hermite in his proof that
$e$ is transcendental. Various irrationality proofs of famous mathematical constants
use Hermite-Pad\'e approximation, even though this may not always be 
obvious.
Proving irrationality can be done by constructing good rational approximants:

\begin{lemma}  \label{lem:1}
Let $x \in \mathbb{R}$. Suppose we can find sequences of integers
$(p_n), (q_n)$ such that
\begin{enumerate}
\item $q_n x - p_n \neq 0$ for all $n \in \mathbb{N}$,
\item $\lim_{n \to \infty} (q_n x-p_n) = 0$.
\end{enumerate}
Then $x$ is irrational.
\end{lemma}

\begin{proof}
Suppose that $x$ is rational. Then $x=p/q$ for some integers $p,q$. We then have
\[   q_n x - p_n = \frac{q_np - p_n q}{q} \]
and since this is not zero for every $n$, we see that $q_np-p_nq \neq 0$ for all $n$.
But since these are integers, this implies that $|q_np-p_nq| \geq 1$ for all $n$.
This shows that $|q_n x - p_n| \geq 1/q$, which is in contradiction with condition 2
in the lemma. Hence we must conclude that $x$ is irrational.
\end{proof}

The construction of the sequences $p_n$ and $q_n$ often uses Pad\'e or Hermite-Pad\'e
approximation for well chosen functions. As an example, consider the two Markov functions
\[   f_1(z) = \int_{0}^1 \frac{dx}{z-x}, \quad f_2(z) =
\int_{-1}^0
\frac{dx}{z-x}, \]
which form an Angelesco system. Some straightforward calculus gives
\[    f_1(i) = -\frac12 \log 2 - \frac{i\pi}{4}, \quad f_2(i) =
\frac12 \log 2 - \frac{i\pi}{4}, \]
hence the sum gives $f_1(i)+f_2(i) = -i\pi/2$. The type II Hermite-Pad\'e approximants
for $f_1$ and $f_2$ will give approximations to $\pi$. Recall that
\begin{eqnarray*}
   P_{n,n}(z) f_1(z) - Q_{n,n;1}(z) & = & \int_0^1
   \frac{P_{n,n}(x)}{z-x} \, dx               \\
   P_{n,n}(z) f_2(z) - Q_{n,n;2}(z) & = & \int_{-1}^0
   \frac{P_{n,n}(x)}{z-x}\, dx .
\end{eqnarray*}
Summing both equations gives
\[  P_{n,n}(z)[f_1(z)+f_2(z)] - [Q_{n,n;1}(z)+Q_{n,n;2}(z)] =
 \int_{-1}^1 \frac{P_{n,n}(x)}{z-x}\, dx . \]
So the fact that we are using a common denominator comes in very handy here.
Then we evaluate these expressions at $z=i$ and hope that
$P_{n,n}(i)$ and $Q_{n,n;1}(i)+Q_{n,n;2}(i)$ are (up to the factor $i$) integers or rational
numbers with simple denominators. Conditions 1 and 2 in Lemma \ref{lem:1} can be checked
by using asymptotic properties of Hermite-Pad\'e approximation. For this particular case
the type II multiple orthogonal polynomials are given by a Rodrigues formula
\[   P_{n,n}(x) = \frac{d^n}{dx^n} \left( x^n(1-x^2)^n \right), \]
and these polynomials are known as \dword{Legendre-Angelesco polynomials}. They have been
studied in detail by Kalyagin \cite{kal} (see also \cite{wva1}). The Rodrigues formula
in fact simplifies the asymptotic analysis, since integration by parts now gives
\[  \int_{-1}^1 \frac{P_{n,n}(x)}{z-x}\, dx = \int_{-1}^1
(-1)^n n! \frac{x^n(1-x^2)^n}{(z-x)^{n+1}}\, dx, \]
which can be handled easily. Some trial and error show that one gets better
results by taking $2n$ instead of $n$, and by differentiating $n$ times:
\begin{multline}  \label{eq:pi}
 \frac{d^n}{dz^n}  \left( P_{2n,2n}(z)[f_1(z)+f_2(z)] -
[Q_{2n,2n;1}(z)+Q_{2n,2n;2}(z)] \right)_{z=i} \\ = (3n)!(-i)^{n+1}
  \int_{-1}^1 \frac{x^{2n}(1-x^2)^{2n}}{(1+ix)^{3n+1}}\, dx.
\end{multline}
This gives rational approximants to $\pi$ of the form
\[  \pi = \frac{b_n}{a_n c_n} + \frac{K_n}{a_n}, \]
where $a_n,b_n,c_n$ are explicitly known integers and $K_n$ is the integral on the
right hand side of (\ref{eq:pi}). The rational approximants show that $\pi$ is irrational
(which was shown already in 1773 by Lambert), and they even show that you can't approximate
$\pi$ by rationals at order greater than $23.271$ (Beukers \cite{Beuk}), i.e.,
\[   \left| \pi - \frac{p}{q} \right| < \frac{1}{q^r}, \]
with $r > 23.271$ only has a finite number of solutions $(p,q)$, where $p$ and $q$ are relatively prime
integers.
This upper bound for the order of approximation
can be reduced to $8.02$ (Hata \cite{Hata}) by considering Markov functions $f_1$ and $f_3$, with
\[   f_3(z) := \int_{-i}^0 \frac{dx}{z-x}. \]
This $f_3$ is now over a complex interval, and then Theorem  \ref{thm:ang}
concerning the location of the zeros no longer holds, and the asymptotic behavior 
must be handled by another method.

One can also use Hermite-Pad\'e approximants to prove transcendence. Then one uses the following
lemma, which extends Lemma \ref{lem:1} from irrational numbers to non-algebraic numbers.
\begin{lemma}  \label{lem:2}
Let $x \in \mathbb{R}$. Suppose that for every integer $m \in \mathbb{N}$ and for all integers
$a_0,a_1,\ldots,a_m \in \mathbb{Z}$ we can find integers $p_{0,n},p_{1,n},\ldots,p_{m,n}$
such that
\begin{enumerate}
\item $\sum_{k=0}^m a_k p_{k,n} \neq 0$ for all $n \in \mathbb{N}$,
\item $\lim_{n \to \infty} (p_{0,n} x^k - p_{k,n}) = 0$ for $k=1,2,\ldots,m$.
\end{enumerate}
Then $x$ is transcendental.
\end{lemma}

\begin{proof}
Suppose that $x$ is algebraic. Then there exists an integer $m$ and integers $a_0,\ldots,a_m$
such that $\sum_{k=0}^m a_kx^k = 0$. But then
\[    \sum_{k=0}^m a_k (p_{0,n}x^k - p_{k,n}) = -\sum_{k=0}^m a_k
p_{k,n} . \]
The right hand side is an integer different from zero, hence
\[   \left| \sum_{k=0}^m a_k( p_{0,n}x^k - p_{k,n}) \right| \geq
1, \]
for all $n \in \mathbb N$. But this contradicts condition 2 of the lemma. Hence
we must conclude that $x$ is not algebraic.
\end{proof}

If we use type II Hermite-Pad\'e approximation to $(e^{\lambda_1x}, e^{\lambda_2x},\ldots, e^{\lambda_rx})$
near $x=0$, then this will give the transcendence of $e$. For Hermite-Pad\'e approximation
near $x=0$ we can use two multi-indices $\vec{n}=(n_1,n_2,\ldots,n_r)$ and $\vec{m} = (m_1,m_2,\ldots,m_r)$.
These Hermite-Pad\'e approximants are known explicitly when 
$m_j+n_j=N+|\vec{n}|$ for $1 \leq j \leq r$, where $N$ is an integer. If
we define the polynomial
\[   T(x) := x^N(x-\lambda_1)^{n_1}(x-\lambda_2)^{n_2}\cdots
(x-\lambda_r)^{n_r}, \]
then $T$ has degree $N+|\vec{n}|$. The expression
\[  P_{\vec{n}}(z) = z^{|\vec{n}|+N+1} \int_0^\infty T(x)
e^{-zx}\, dx \]
gives a polynomial of degree $|\vec{n}|$, and
\[    Q_{\vec{m};j}(z) = z^{|\vec{n}|+N+1} \int_0^\infty
T(x+\lambda_j) e^{-zx}\, dx \]
gives a polynomial of degree $|\vec{n}| + N -n_j=m_j$. One easily verifies
that
\[   P_{\vec{n}}(z) e^{\lambda_j z} - Q_{\vec{m};j}(z) =
 e^{\lambda_j z} z^{|\vec{n}|+N+1} \int_0^{\lambda_j} T(x)
 e^{-zx}\, dx = \O(z^{n_j+m_j+1}),  \]
as $z \to 0$,  which are the interpolation conditions for type II Hermite-Pad\'e
 approximation near the origin for the two multi-indices $(\vec{n},\vec{m})$.

 For proving the transcendence of $e$, we take $\lambda_j=j$, $z=1$ and
 for a prime $p > r$, which is not a divisor of $a_0$, we take $N=p-1$ and $n_j=p$ $(j=1,\ldots,r)$.
 Then some elementary calculus shows that $p_0=P_{\vec{n}}(1)/(p-1)!$
 is an integer which is not divisible by $p$ and each
 $p_j = Q_{\vec{m};j}(1)/(p-1)!$ is an integer divisible by $p$. Therefore
 $\sum_{j=0}^r a_j p_j$ is not divisible by $p$ and hence
 condition 1 of Lemma \ref{lem:2} is satisfied.
 Furthermore
 \[   p_0 e^j - p_j = \frac{e^j}{(p-1)!} \int_0^j T(x) e^{-x}\,
 dx, \]
 and the simple estimate $|T(x)|\leq j^{(r+1)p-1}$ on $[0,j]$,
 shows that this converges to $0$ for every $j$ when the prime $p$
 tends to infinity (luckily Euclides showed that there are
 infinitely many primes). So condition 2 of Lemma \ref{lem:2} is
 also satisfied and we conclude that $e$ is transcendental
 (Hermite, 1874).

\subsection{Other applications}
Recently a number of applications came up in other areas of mathematics and theoretical physics.
There are interesting connections with random matrix theory, where multiple orthogonal polynomials
(in particular multiple Hermite polynomials) appear when one investigates random matrices
with an external source \cite{BK1,ABK}. Multiple Laguerre polynomials appear for the Wishart
ensemble of random matrices \cite{BK}. Multiple Jacobi polynomials (the Jacobi-Pi\~neiro polynomials)
were used to obtain a counterexample to the Bethe Ansatz Conjecture for the Gaudin model \cite{MV}.
More details on multiple orthogonal polynomials (recursion relation, specific examples, etc.)
can be found in \cite[Chapter 23]{ism}.

\subsection*{Acknowledgments} Research supported by research grant OT/04/21 of Katholieke
Universiteit Leuven, research project G.0455.04 of FWO-Vlaanderen and INTAS research network
03-51-6637.

{\obeylines
Walter Van Assche
Katholieke Universiteit Leuven
Department of Mathematics
Celestijnenlaan 200B
B-3001 Leuven
BELGIUM
\tt walter@wis.kuleuven.be
http://www.wis.kuleuven.be/analyse/walter/}
\endddoc